\documentclass[12pt]{article}
\usepackage {amsfonts}
\usepackage[russian,english]{babel}

\newtheorem{Th}{Theorem}
\newtheorem{Le}{Lemma}
\newtheorem{De}{Definition}

\def\R{{\Bbb R}}
\def\C{{\Bbb C}}

\def\Z{{\Bbb Z}}
\def\Q{{\Bbb Q}}
\def\e{\varepsilon}
\def\ph{\varphi}
\def\O{\Omega}
\def\p{\varphi}
\def\<{\langle }
\def\>{\rangle }
\def\H{{\cal H}}
\def\D{{\cal D}}
\def\bp{\bar\partial}
\def\p{\partial}
\def\eq{eqnarray}
\def\l{\lambda}
\def\a{\alpha}
\def\b{\beta}
\def\g{\gamma}
\def\dist{{\rm dist}}

\begin{document}

\title{
    \bf Zeroes of holomorphic functions with almost--periodic modulus}

\author{Favorov S.Ju.}

\date{}

\maketitle

\begin{abstract}

We give necessary and sufficient conditions for a divisor in a
tube domain to be the divisor of a holomorphic function with
almost--periodic modulus.

\medskip


\medskip


\end{abstract}

\medskip

Zero distribution for various classes of holomorphic
almost--periodic functions in a strip was studied by many authors
(cf. \cite{B}, \cite{FRR1}, \cite{J}, \cite{JT}, \cite{KL},
\cite{L}, \cite{T1}). The notion of almost--periodic discrete set
appeared in \cite{KL} and \cite{T1} in connection with these
investigations. Its generalization to several complex variables
was the notion of almost--periodic divisor, introduced by
L.\,I.\,Ronkin (cf. \cite{R1}) and studied in his works and works
of his disciples (cf. \cite{FRR2}, \cite{FRR3}, \cite{R4}). But
these notions are not sufficient for a complete description of
zero sets of holomorphic almost--periodic functions (cf.
\cite{T2}): in addition, one needs some topological
characteristic, namely, Chern class of the special (generated by
an almost--periodic set or a divisor) line bundle over Bohr's
compact set (cf. \cite{F1}, \cite{F2}). On the other hand, the
class of zero sets of holomorphic functions with almost--periodic
modulus in a strip is just the class of almost--periodic discrete
sets (cf. \cite{FRR1}). That's why it is natural to obtain a
description of zeroes of holomorphic functions with the
almost--periodic modulus for several complex variables without
using topological terms. This problem is just solved in our paper.
\medskip

By $T_S$ denote a tube set $\{z=x+iy:\,x\in\R^m, y\in S\}$, where
the base $S$ is a subset of $\R^m$.

\begin{De}. A continuous function $f$ on $T_S$ is called almost--periodic,
if for each sequence $\{f(z+h_n)\}_{h_n\in\R^m}$ of shifts there
exists a uniformly convergent on $T_S$ subsequence.
\end{De}

In particular, for $S=\{0\}$ we obtain the definition of an
almost--periodic function on $\R^m$. \footnote {This definition is
equivalent to another one that makes use of the notion of an
$\e$-almost period; for $m=1$ see, for example, \cite{Lt}, the
extension to $m>1$ is trivial.}

It follows easily that any almost--periodic function on a tube set
with a compact base is bounded.

\begin{De}. Let $\O$ be a domain in $\R^m$. A continuous function $f$ on
$T_\O$ is called almost--periodic, if its restriction to every
tube set $T_K$ with compact base $K\subset\O$ is an
almost--periodic function on  $T_K$.
\end{De}

\begin{De}\label{def} (cf. \cite{R1}, for distributions from $\D'(\R)$
cf. also \cite{S}). A distribution $F(z)\in\D'(T_\O)$ is called
almost--periodic, if for any test--function $\ph(z)\in\D(T_\O)$
the function $\<F(z),\ph(z-t)\>$ is an almost--periodic function
in $t\in\R^m$.
\end{De}

The next assertion is valid.
\begin{Th}\label{equ} (cf. \cite{R1}).
A distribution $F\in\D'(T_\O)$ is almost periodic if and only if
for each sequence $\{h^n\}\in\R^m$ there exists a subsequence
$\{\tilde h^n\}$ such that the sequence of the distributions
$F(z+\tilde h^n)$ converges uniformly on the sets
$\{\kappa(z-t):\,t\in\R^m,\,\kappa\in{\cal K}\}$, where ${\cal K}$
is any compact family in $\D(T_\O)$.
\end{Th}

\begin{De}\label{mean} (cf. \cite{R1})
The mean value (in the variable $x\in\R^m$) of an almost--periodic
distribution $F$ is the distribution $c_F(y)\otimes dx$ with
$c_F(y)\in\D'(\O)$ and the Lebesgue measure $dx$ on $\R^m$,
defined for a test--functions $\ph\in\D(T_\O)$ by the equality
$$
\<c_F(y)\otimes
dx,\ph(z)\>=\lim_{N\to\infty}(2N)^{-m}\int_{\max_j|t_j|<N}
\<F(z),\ph(z-t)\>dt,
$$
where $t=(t_1,\dots,t_m)\in\R^m$.
\end{De}

Note that if $F(z)$ is an almost--periodic function on $T_\O$,
then $c_F(y)$ is a continuous function on $\O$. Further, if $F(z)$
is an almost--periodic complex measure on $T_\O$, then $c_F(y)$ is
a complex measure on $\O$ as well, and $c_F(y)\otimes dx$ is the
weak limit of the measures $F(tx+iy)dx\,dy$ as $|t|\to\infty$ (cf.
\cite{R1}).

\medskip

By $\H(G)$ denote the space of holomorphic functions on the domain
$G\subset\C^m$ with respect to the topology of the uniform
convergence on compact subsets of $G$.

The following assertion is true.

\begin{Th}\label{tg} (cf. \cite{R1}).
If a function $f\in\H(T_\O)$ is almost--periodic, then $\log |f|$
is an almost--periodic distribution on $T_\O$.
\end{Th}

The main part of the proof of this theorem is the following lemma.

\begin{Le}\label{con}(cf. \cite{R1}).
If $f_n\in\H(G),\, n=1,2,\dots,$ and $f_n\to f\not\equiv 0$ in the
space $\H(G)$, then $\log |f_n|\to\log|f|$ in the space $\D'(G)$.
\end{Le}

Now let
\begin{\eq}\label{cur}
(i/\pi)\p\bp\log|f|=(2/\pi)\sum_{j,k=1}^m{\p^2\log|f|\over\p
z_j\p\bar z_k}(i/2)dz_j\wedge d\bar z_k.
\end{\eq}
be the current of integration over the divisor $d_f$ of the
function $f(z)\in\H(G)$, $z=(z_1,\dots,z_m)$. In the case $m=1$
this current corresponds to the discrete measure with integer
masses equal to the multiplicities of the zeroes of the function
$f$.

Note that all the coefficients of the current (\ref{cur}) are
complex measures on $G$, and the "diagonal" coefficients
$\p^2\log|f|\over\p z_j\p\bar z_j$ are positive measures (cf.
\cite{LG}).

Clearly, the differentiation keeps the almost periodicity of
distributions. Therefore, it follows from Theorem \ref{tg} that
all the coefficients of the current (\ref{cur}) are
almost--periodic distributions for any holomorphic
almost--periodic function on $T_\O$. If we replace $f$ by another
holomorphic function on $T_\O$ with the same divisor, then the
coefficients of the current (\ref{cur}) do not change. Hence an
almost--periodicity of all the coefficients does not imply almost
periodicity of the function $f$ itself.

\begin{De}
(cf. \cite{FRR2}, \cite{FRR3}). The divisor $d_f$ of a function
$f\in\H(T_\O)$ is called almost--periodic, if all the coefficients
of the current (\ref{cur}) are almost--periodic distributions.
\end{De}

Note that in \cite{R1} a divisor $d_f$ was called
almost--periodic, if the measure $\sum_{j=1}^m{\p^2\log|f|\over\p
z_j\p\bar z_j}$, was almost--periodic on $T_\O$. But that
definition is equivalent to the given above (cf. \cite{FRR3}).

\smallskip

There exist almost periodic divisors which cannot be generated by
holomorphic almost periodic functions. For example, let $g(w)$ be
an entire function on $\C$ with simple zeroes at the points of the
standard integer--valued lattice, and let $d[\l,\mu]$,
$\l,\mu\in\R^m$ be the divisor of the function
$g(\<z,\,\l\>+i\<z,\,\mu\>)$. This divisor is periodic for
 vectors $\l,\,\mu$ that are linearly dependent over $\Q$ or linearly
 independent over $\R$ (with the periods ${|\mu|^2\l-\<\l,\,\mu\>\mu\over
 |\l|^2|\mu|^2-\<\l,\,\mu\>^2}$ and
${|\l|^2\mu-\<\l,\,\mu\>\l\over |\l|^2|\mu|^2-\<\l,\,\mu\>^2}$).
Then $d[\l,\mu]$ is almost periodic for $\l,\,\mu$ linearly
independent over $\Q$ and linearly dependent over $\R$ (for $m=1$
cf. \cite{T2}; since a real linearly transform in $\C^m$ keeps
almost--periodicity, the case $m>1$ follows as well). Besides, the
divisor $d[\l,\mu]$ for any linearly independent over $\Q$ vectors
$\l,\,\mu$ is the divisor of no holomorphic almost periodic
function (in the case $m=1$, i.e., irrational $\l/\mu$
cf.\cite{T2}, for $m>1$ cf. \cite{R4}). A complete description of
the divisors of holomorphic almost--periodic functions is
contained in the following theorem.

\begin{Th}\label{coh}(for $m=1$ cf. \cite{F1}, for $m>1$ cf. \cite{F2}).
A holomorphic bundle over Bohr's compactification $K_B$ of the
space $\R^m$ is assigned to each almost--periodic divisor $d$ on a
tube domain $T_\O$ with convex base $\O$ such that:

the map $d\mapsto c(d)$, $c(d)$ being the first Chern class of
this bundle, is a homomorphism of the semigroup of positive
almost--periodic divisors on $T_\O$ to the cohomology group
$H^2(K_B,\Z)$, the kernel of this homomorphism is just the set of
all divisors of holomorphic almost--periodic functions on $T_\O$,

a finite family $\l^j,\mu^j\in\R^m$ corresponds to each cohomology
class $c(d)$ such that $c(d)=\sum_jc(d[\l^j,\mu^j])$,

the mapping $W:\,(\l,\mu)\mapsto c(d[\l,\mu])$ is skew-symmetric
and additive in variables $\l,\mu\in\R^m$.
\end{Th}

A description of zeroes for holomorphic functions of one variable
with the almost--periodic modulus is given in the following
theorem.

\begin{Th}\label{tr}
(cf. \cite{FRR1}; for divisors $d[\l,\mu],\,\l,\mu\in\R$ cf.
\cite{T2}). A divisor $d$ on a strip is the divisor of some
holomorphic function on the strip with almost--periodic modulus if
and only if $d$ is almost--periodic.
\end{Th}

Now consider the multidimensional case again. Note that for an
almost--periodic divisor $d$ on  $T_\O$ all the coefficients of
the current (\ref{cur}) have mean values in $x$. The imaginary
parts of these mean values, i.e., the mean values of the real
measures $(2/\pi)\Im{\p^2\log|f|\over\p z_j\p\bar z_k}$ have the
form $a_{j,k}dy\otimes dx$, $a_{j,k}\in\R$ (cf. \cite{FRR3}). By
$A(d)$ denote the matrix with the entries $a_{j,k}$. In the case
$d=d_f$ for an almost--periodic function $f\in\H(T_\O)$ we have
$A(d)=0$ (cf. \cite{R2}).

\begin{Th}\label{main}.
A divisor $d$ on a tube domain $T_\O$ with convex base $\O$ is the
divisor of a holomorphic function with almost--periodic modulus if
and only if divisor $d$ is almost--periodic, and the
skew-symmetric matrix $A(d)$ is zero.
\end{Th}

To prove this theorem we need the following improvement of Theorem
\ref{tg}.

\begin{Th}\label{tm}.
A function $f\in\H(T_\O),\,f\not\equiv 0$, has almost--periodic
modulus if and only if the distribution $\log|f|\in\D'(T_\O)$ is
almost--periodic.
\end{Th}

{\bf Proof of Theorem \ref{tm}}. Let $|f(z)|$ be an
almost--periodic function on $T_\O$, and let $\{h^n\}$ be an
arbitrary sequence from $\R^m$. In order to check that $\log|f|$
is an almost--periodic distribution on $T_\O$, we will prove that
for any continuous function $\ph$ with compact support in $T_\O$,
the sequence of functions
\begin{\eq}\label{c}
 \psi_n(t)=\int\log|f(z+h^n)|\ph(z-t)dx dy
\end{\eq}
contains a convergent, uniformly on $\R^m$, subsequence. We will
prove this assertion by contradiction.

First, since the function $|f(z)|$ is uniformly bounded on $T_K$
for every compact set $K\subset\O$, we may assume that the
sequence of the functions $\{f(z+h^n)\}$ converges to some
function $g(z)$ in the space $\H(T_\O)$. Further, since the
function $|f(z)|$ is almost--periodic on $T_\O$, we may assume
that the sequence of the functions $\{|f(z+h^n)|\}$ converges to
some function $\Phi(z)\not\equiv 0$ uniformly on each $T_K$. If
the sequence (\ref{c}) does not converge uniformly on $\R^m$, then
for some $\delta>0$ and some subsequence of $n$ there exist
$t^n\in \R^m$ with the property
\begin{\eq}\label{ge}
|\psi_n(t^n)-\int\log|g(z)|\ph(z-t^n)dx dy|\ge\delta.
\end{\eq}
The function $|g(z)|\equiv\Phi(z)$ is almost--periodic on $T_\O$,
hence we may assume that the same subsequence of the functions
$\{|g(z+t^n)|\}$ converges uniformly on each $T_K$ to some
function $\Psi(z)\not\equiv 0$. Since the sequence of the
functions $\{|f(z+h^n+t)|\}$ converges uniformly in $t\in\R^m$ and
$z\in T_K$ to the function $|g(z+t)|$, we see that the subsequence
$\{|f(z+h^n+t^n)|\}$ converges to $\Psi(z)$ uniformly on $T_K$.
Also, the subsequences of the functions $\{f(z+h^n+t^n)\}$ and
$\{g(z-t^n)\}$ are bounded uniformly on compact subsets of $T_\O$,
therefore passing to a subsequence again, we get that
$f(z+h^n+t^n)\to H_1(z)$ and $g(z+t^n)\to H_2(z)$ in the space
$\H(T_\O)$, and $|H_1(z)|=\Psi(z)=|H_2(z)|$. Using Lemma
\ref{con},  we obtain that the corresponding subsequences of the
functions $\{\log|f(z+h^n+t^n)|\}$ and $\{\log|g(z+t^n)|\}$
converge, in the sense of distributions, to the same function
$\log\Psi(z)$. The last assertion contradicts (\ref{ge}).

On the other hand, let $\log|f(z)|$ be an almost--periodic
distribution on $T_\O$, and let $\ph_\e(z)$ be a nonnegative,
depending on $|z|$ smooth function such that $\ph(z)=0$ for
$|z|>\e$ and $\int_{\C^m}\ph_\e(z)dx\,dy=1$. Evidently, the family
of functions $\{\ph_\e(z+iy)\}_{|y|\le C}$ is a compact set in the
space $\D(\C^m)$ for every $C<\infty$. Let $K$ be a compact set in
$\O$ and $\e<\dist\{K,\p\O\}$. Now Theorem \ref{equ} implies that
the convolution $(\log|f|\ast\ph_\e)(z)$ is an almost--periodic
function on $T_K$. Hence this convolution is bounded on $T_K$, and
the inequality $\log|f(z)|\le(\log|f|\ast\ph_\e)(z)$ shows that
$|f(z)|$ is bounded on $T_K$ as well.

Suppose that $|f|$ is not an almost--periodic function on $T_\O$.
Then there exists a sequence of functions
$\{|f(z+h^n)|\},\,h^n\in\R^m,$ such that every its subsequence
does not converge uniformly on $T_{K'}$ for some compact set
$K'\subset\O$. Without loss of generality it can be assumed that
the sequence of functions $\{f(z+h^n)\}$ converges in the space
$\H(T_\O)$ to some function $g(z)$. It is clear that $g(z)$ is
bounded on $T_K$ for every compact set $K\subset\O$. Further, by
Lemma \ref{con} we get $\log|f(z+h^n)|\to\log|g(z)|$ in the sense
of distributions. Using Theorem \ref{equ} and passing to a
subsequence, we obtain
\begin{\eq}\label{ge2}
\int(\log|f(z+h^n)|-\log|g(z)|)\ph_\e(z-t-is)dxdy\to 0.
\end{\eq}
uniformly in $t\in\R^m$ and $s\in K'$. On the other hand, for some
$\delta>0$ and some subsequence of $n$ there exist points
$z^n=x^n+iy^n\in T_K'$ such that
\begin{\eq}\label{ge1}
||f(h^n+x^n+iy^n)|-|g(x^n+iy^n)||\ge\delta.
\end{\eq}
Passing to a subsequence if necessary, we may assume that $y^n\to
y^0\in K'$, and the sequences of the functions
$\{f(z+h^n+z^n-iy^0)\}$ and $\{g(z+z^n-iy^0)\}$ converge in the
space $\H(T_\O)$ to functions $H_1(z)$ and $H_2(z)$, respectively.
Then Lemma \ref{con} implies that
 $\log|f(z+h^n+z^n-iy^0)|\to\log|H_1(z)|$ and
$\log|g(z+z^n-iy^0)|\to\log|H_2(z)|$ in the space $\D'(T_\O)$.
Taking into account (\ref{ge2}), we obtain
$$
\int\log|H_1(z)|\ph_\e(z-iy^0)dx\,dy=\int\log|H_2(z)|\ph_\e(z-iy^0)dxdy.
$$
Since $\e$ is arbitrary small, we get $|H_1(iy^0)|=|H_2(iy^0)|$.
At the same time, by (\ref{ge1}) we have $|H_1(iy^0)|\neq
|H_2(iy^0)|$. This contradiction proves Theorem \ref{tm}.

\medskip

{\bf Proof of the necessity in Theorem \ref{main}}. It follows
from Theorem \ref{tm} that every function  $f\in\H(T_\O)$ with
almost--periodic modulus has an almost--periodic divisor. Further,
the mean value $c_{\log|f|}(y)\otimes dx$ of the function
$\log|f|$ is the weak limit of the measures
$\log|f(tx+iy)|dx\otimes dy$ as $|t|\to\infty$ in the space
$\D'(T_\O)$, therefore for all $j\neq k$ the mean values of the
distributions
$$
\Im{\p^2\log|f|\over\p z_j\p\bar z_k}= {1\over
4}\left({\p^2\over\p x_j\p y_k}-{\p^2\over\p x_k\p y_j}\right)
\log|f|
$$
equal
\begin{\eq*}
\lim_{|t|\to\infty} {1\over 4}\left({\p^2\over\p x_j\p
y_k}-{\p^2\over\p x_k\p y_j}\right)\log|f(tx+iy)|dx\otimes dy=\\
{1\over 4}\left({\p^2\over\p x_j\p y_k}-{\p^2\over\p x_k\p
y_j}\right)c_{\log|f|}(y)\otimes dx=0.
\end{\eq*}
The necessity of the conditions in Theorem \ref{main} is proved.

\medskip

The proof of the sufficiency makes use of the following lemmas. As
above, $d[\l,\mu]$, $\l,\mu\in\R^m$ is the divisor of the function
$g(\<z,\,\l\>+i\<z,\,\mu\>)$, where $g(w)$ is an entire function
on $\C$ with simple zeroes at the points of the standard
integer--valued lattice.

\begin{Le}\label{mod}.
The divisor $d[\l,\mu]$ with $\l=t\mu,\,\l\in\R^m,\,t\in\R$, is
the divisor of an entire function on $\C^m$ with almost--periodic
modulus.
\end{Le}

{\bf Proof of Lemma \ref{mod}}. After a suitable regular real
linear transform we obtain the case $\mu=(1,0,\dots,0)$, i.e., the
case of a divisor depending only on one coordinate, therefore the
assertion of our lemma is a consequence of Theorem \ref{tr}.

\medskip

 Further, let $e^1,\dots,e^m$ be the coordinate  vectors
in $\C^m$.

\begin{Le}\label{dis}.
The entries $a_{j,k}$ of the matrix $A_0=A(d[e^1,e^2])$ vanish for
$(j,\,k)\neq(1,\,2)$ or $(2,\,1)$, and $a_{1,2}=-1,\,a_{2,1}=1$.
\end{Le}

{\bf Proof of Lemma \ref{dis}}. The divisor of the function
$g(z_1+iz_2)$ does not depend on variables $z_j$ with $j>2$, hence
the distributions $\Im{\p^2\log|g(z_1+iz_2)|\over\p z_j\p\bar
z_k}$ vanish for $(j,\,k)\neq(1,\,2)$ or $(2,\,1)$.

Consider the expression
\begin{\eq}\label{d}
(L_z\log|g(z_1+iz_2)|,\,\ph(z_1+t_1,z_2+t_2)),\quad
(t_1,\,t_2)\in\R^2,
\end{\eq}
for $L_z={2\over\pi}\Im{\p^2\over\p\bar z_1\p z_2}$ and a function
$\ph(z)\ge 0$ from the space $\D(\C^2)$. In the coordinates
$\zeta_1=z_1+iz_2,\;\zeta_2=z_1-iz_2$, it has a form
$$
{1\over 4}(\tilde L_\zeta\log|g(\zeta_1)|,
\,\ph((\zeta_1+\zeta_2)/2+t_1,(\zeta_1 -\zeta_2)/2i+t_2))
$$
with
$$
\tilde
L_\zeta={2\over\pi}\Re\left({\p^2\over\p\zeta_1\p\bar\zeta_1}-
{\p^2\over\p\zeta_2\p\bar\zeta_1}+{\p^2\over\p\zeta_1\p\bar\zeta_2}-
{\p^2\over\p\zeta_2\p\bar\zeta_2}\right).
$$
Using the definition of $g$ and properties of the Laplace
operator, we get
\begin{\eq*}
\tilde L_\zeta\log|g(\zeta_1)|=
{2\over\pi}{\p^2\over\p\zeta_1\p\bar\zeta_1}\log|g(\zeta_1)|=
\sum_{q_1,q_2\in\Z}\delta(\zeta_1-q_1-iq_2)\otimes d\xi d\eta,
\end{\eq*}
where $\delta$ is the Dirac function on the plane,
$\xi=\Re\zeta_2, \eta=\Im\zeta_2$. Therefore, (\ref{d}) is equal
to
$$ {1\over
4}\sum_{q_1,q_2\in\Z}\int_\C\ph(t_1+(q_1+iq_2+\xi+i\eta)/2,
t_2+(q_1+iq_2-\xi-i\eta)/2i)d\xi d\eta.
$$
Substituting $\xi-q_1,\;\eta+q_2$ for $u,\,v$, respectively, we
get
\begin{\eq}\label{di}
{1\over 4}\sum_{q_1,q_2\in\Z}\int_\C\ph(u/2+iv/2+t_1+q_1,-v/2+iu/2
+t_2+q_2)du\,dv.
\end{\eq}
Since the divisor $d_{e^1,e^2}$ has period $1$ in each variable,
we see that the mean value of (\ref{d}) is the integral of
(\ref{di}) over the square $0\le t_1\le 1,\;0\le t_2\le 1$. Then
it is equal to the integral
$$
{1\over 4}\int_{\R^4}\ph(u/2+iv/2+t_1,-v/2+iu/2
+t_2)du\,dv\,dt_1\,dt_2.
$$
Finally, substituting
$u/2+t_1=x_1,\;v/2=y_1,\;t_2-v/2=x_2,\;u/2=y_2$, we obtain the
equality
\begin{\eq*}
\int_0^1\int_0^1(L_z\log|g(z_1+iz_2)|,\,\ph(z_1+t_1,z_2+t_2))dt_1\,dt_2=\\
\int_{\R^4}\ph(x_1+iy_1,x_2+iy_2)dx_1\,dy_1\,dx_2\,dy_2,
\end{\eq*}
hence the mean value of the distribution $L_z\log|g(z_1+iz_2)|$ is
the Lebesgue measure in $\C^2$. The lemma is proved.

\medskip

By $(\l,\mu)$ denote the matrix product $(\l_j\mu_k)_{j,k=1}^m$ of
the vectors
$\l=(\l_1,\dots,\l_m),\;\mu=(\mu_1,\dots,\mu_m)\in\R^m$.

\begin{Le}\label{div}.
For any $\l,\,\mu\in\R^m$, the matrix $A(d[\l,\mu])$ equals the
difference $(\mu,\l)-(\l,\mu)$.
\end{Le}

{\bf Proof of Lemma \ref{div}}. If $\l,\,\mu$ are linearly
dependent over $\R$, then $(\mu,\l)-(\l,\mu)=0$. On the other
hand, it follows from Lemma \ref{mod} that in this case the
divisor $d[\l,\mu]$ is the divisor of some holomorphic in $\C^m$
function with almost--periodic modulus. Using the proved part of
Theorem \ref{main}, we have $A(d[\l,\mu])=0$.

Let $\l,\,\mu$ be linearly independent over $\R$. The divisor
$d[\l,\mu]$ is the divisor $d[e^1,e^2]$ in the coordinates
$\zeta=Bz$ for some real nondegenerate matrix $B$ with the first
and second rows $\l$ and $\mu$, respectively.  Note that the
matrix $A(d)$ is the matrix of the mean values for the matrix
${1\over 2i}(D(z)-\bar D(z))$, where
$$
D(z)=\left({\p^2\log|g(\<z,\l\>+i\<z,\mu\>)|\over\p z_j\p\bar
z_k}\right),
$$
$\bar D$ being the  matrix with all the entries complex conjugated
to the corresponding entries of the matrix $D$. Therefore
$D(z)=B'D(\zeta)B$, $B'$ being the transpose matrix to $B$, and
$A(d[\l,\mu])=B'A_0 B$ for the matrix $A_0$ from the previous
lemma. This completes the proof of Lemma \ref{div}.

\medskip

\begin{Le}\label{w1}.
If numbers $\a_j,\b_j\in\R,\;j=1,\dots,n$, satisfy the condition
$\sum_1^n\a_j\b_j=0$, then for some
$\g_k\in\R,\;\nu^k\in\R^m,\;k=1,\dots,N$, we get
\begin{\eq}\label{eq}
\sum_1^nW(\a_je^1,\,\b_je^2)=\sum_1^N W(\g_k\nu^k,\,\nu^k),
\end{\eq}
$W$ being the mapping from Theorem \ref{coh}.
\end{Le}

{\bf Proof of Lemma \ref{w1}}. The case $n=1$ means that the
left-hand side of (\ref{eq}) vanishes. For $n>1$ we have
$$
W(\a_{n-1}e^1,\,\b_{n-1}e^2)+W(\a_ne^1,\,\b_ne^2)=
W(\a_{n-1}e^1,\,\a_ne^1)+W(\b_ne^2,\,\b_n\a_n/\a_{n-1}e^2)
$$
$$
+W(\a_ne^1+\a_n\b_n/\a_{n-1}e^2,\,\a_{n-1}e^1)+\b_ne^2)+
W(\a_{n-1}e^1,\,(\b_{n-1}+\b_n\a_n/\a_{n-1})e^2).
$$
The first three terms of the right-hand side have the form
$W(\g\nu,\,\nu),\,\g\in\R, \nu\in\R^m$. Subtracting these terms
from the left-hand side of (\ref{eq}), we get
$$
\sum_1^{n-2}W(\a_je^1,\,\b_je^2)+
W(\a_{n-1}e^1,\,(\b_{n-1}+\b_n\a_n/\a_{n-1})e^2).
$$
Hence the lemma can be proved by induction over $n$.

\begin{Le}\label{w2}.
Let vectors $\l^j,\mu^j\in\R^m,\;j=1,\dots,n$ be such that the
matrix $\sum_1^n(\l^j,\mu^j)$ is symmetric. Then
\begin{\eq}\label{eq1}
\sum_1^n W(\l^j,\,\mu^j)=\sum_1^N W(\g_k\nu^k,\,\nu^k)
\end{\eq}
for some $\g_k\in\R,\;\nu^k\in\R^m,\;k=1,\dots,N$.
\end{Le}

{\bf Proof of Lemma \ref{w2}}. The vectors $\l^j,\mu^j$ are linear
combinations of the vectors $e^1,\dots,e^m$, therefore the
left-hand side of (\ref{eq1}) has the form
\begin{\eq}\label{bas}
 \sum_{1\le p,q\le m}
\left(\sum_{j=1}^{M(p,q)}W(\a_{j,p}e^p,\,\b_{j,q}e^q)\right)
\end{\eq}
with $\a_{j,p},\;\b_{j,q}\in\R$. The mapping $W$ is
skew-symmetric, hence we may assume that  all the terms in
(\ref{bas}) vanish for $p>q$, and the entries of the corresponding
matrix
$\left(\sum_{j=1}^{M(p,q)}\a_{j,p}\b_{j,q}\right)_{p,q=1}^m$
vanish for all $p>q$. Since this matrix coincides with the
symmetric matrix $\sum_1^n(\l^j,\mu^j)$, we see that
$\sum_{j=1}^{M(p,q)}\a_{j,p}\b_{j,q}=0$ for $p<q$ as well. Now it
follows from Lemma \ref{w1} that for $p<q$ the sum
$$
\sum_{j=1}^{M(p,q)}W(\a_{j,p}e^p,\,\b_{j,q}e^q)
$$
has the form of the right-hand side of (\ref{eq1}). The terms of
(\ref{bas}) with $p=q$ have already the form $W(\g\nu,\,\nu)$. The
lemma is proved.

\medskip

{\bf Proof of the sufficiency in Theorem \ref{main}}. Let $d$ be a
divisor in $T_\O$ such that $A(d)=0$. It follows from Theorem
\ref{coh} that there exist $\l^j,\,\mu^j\in\R^m,\,j=1,\dots,n$,
such that the sum $d+\sum_j d[\l^j,\mu^j]$ is the divisor of a
holomorphic almost--periodic function. Now, by \cite{R2},
 $A(d+\sum_1^n d[\l^j,\mu^j])=0$. Since the mapping $d\mapsto A(d)$
is a homomorphism, we get $\sum_1^n
(\l^j,\,\mu^j)-(\mu^j,\,\l^j)=\sum_1^n A(d[\l^j,\mu^j])=0$, i.e.,
the matrix $\sum_1^n (\l^j,\,\mu^j)$ is symmetric. Using Lemma
\ref{w2}, we get (\ref{eq1}) for some
$\g_k\in\R,\;\nu^k\in\R^m,\;k=1,\dots,N$. Therefore,
$$
c(d+\sum_1^N d[\g_k\nu^k,\nu^k])=c(d)+\sum_1^N
W(\g_k\nu^k,\,\nu^k)
$$
$$
= c(d)+\sum_1^n W(\l^j,\,\mu^j)=c(d+\sum_1^n d[\l^j,\mu^j])=0.
$$
An application of Theorem \ref{coh} yields that there exists an
almost--periodic function $F\in\H(T_\O)$ with the divisor
$d+\sum_1^N d[\g_k\nu^k,\nu^k]$. Using Lemma \ref{mod}, we can
take functions $f_k\in\H(T_\O)$ with the divisors
$d[\g_k\nu^k,\,\nu^k]$ and almost--periodic modula. The function
$f(z)=F(z)(\prod_1^N f_k(z))^{-1}$ is holomorphic on $T_\O$ and
has the divisor $d$. Then Theorem \ref{tm} implies that the
distributions $\log|F|$ and $\log|f_k|,\;k=1,\dots,N$, are
almost--periodic. Hence the distribution
$\log|f|=\log|F|-\sum_1^N\log|f_k|$ is almost--periodic as well.
Using Theorem \ref{tm} again, we see that the function $|f|$ is
almost--periodic. This completes the proof of Theorem \ref{main}.

\bigskip

Department of Mathematics

Kharkov National University

 Svobody sq.,4, Kharkov 61077

Ukraine

e-mail: favorov@assa.kharkov.ua


\begin{thebibliography}{}

\bibitem{B}
H.\,Bohr, {\it Zur Theorie der Fastperiodischen Funktionen, III
Teil; Dirichletentwicklung Analytischer Funktionen,} Acta math.
{\bf 47} (1926), 237-281.

\bibitem{F1}
Favorov, S.Yu. {\it Zeros of holomorphic almost periodic
functions,} Journal d'Analyse Mathematique 84, (2001), 51-66.


\bibitem{F2}
Favorov, S.Yu., {\it Almost periodic divisors, holomorphic
functions, and holomorphic mappings,} Bull.Sci.Math. 127, (2003),
859-883.

\bibitem{FRR1}
Favorov, S.Yu., Rashkovskii, A.Ju., and Ronkin, L.I., {\it Almost
periodic divisors in a strip,} Journal d'Analyse Math. 74, (1998),
325-345.

\bibitem{FRR2}
Favorov, S.Yu., Rashkovskii, A.Ju., and Ronkin, L.I., {\it Almost
periodic currents and holomorphic chains,} C. R. Acad. Sci. Paris
327, Serie I (1998), 302-307.

\bibitem{FRR3}
Favorov, S.Yu., Rashkovskii, A.Ju., and Ronkin, L.I., {\it Almost
periodic currents, divisors and holomorphic chains in tube
domains,} Israel. Math. Conf. Proc. 15, (2001), 67-88.

\bibitem{J}
B.\,Jessen, {\it \"Uber die Nullstellen einer analitischen
fastperiodischen Funktions, Eine Verallagemeinerung der
Jensenschen Formel},  Math. Ann. {\bf 108} (1933), 485-516.

\bibitem{JT}
B.\,Jessen, H. Tornehave, {\it Mean motions and zeros of almost
periodic functions}, Acta Math. {\bf 77} (1945), 137-279.

\bibitem{KL}
Krejn,M.G., and Levin,B.Ja., {\it On almost periodic functions of
exponential type}, Dokl. AN SSSR {\bf 64} (1949), no. 3, 285-287
(Russian).

\bibitem{L}
Levin,B.Ja., { Distribution of Zeros of Entire Functions.} Transl.
of Math. Monograph, Vol.5, AMS Providence, no. 1, MR 19, 1964,
403p.

\bibitem{LG}
Lelong, P., and Gruman, L., {\it Functions of Several Complex
Variables,} Springer-Verlag. Berlin-Heidelberg, 1986.

\bibitem{Lt}
Levitan,B.M., Almost periodic functions. - M.: Gostehizdat,
 1953, 396 p, (Russian).

\bibitem{R2}
Ronkin,l.I., {\it Jessen's theorem for holomorphic almost periodic
mappings,} Ukrainsk. Mat. Zh. {\bf 42} (1990), 1094-1107
(Russian).

\bibitem{R1}
Ronkin,L.I., {\it Almost periodic distributions and divisors in
tube domains}, Zap. Nauchn. Sem. POMI {\bf 247} (1997), 210-236
(Russian).

\bibitem{R4}
Ronkin, L.I., {\it Almost periodic divisors and the spatial
extension of almost periodic measures}, Linear Topological Spaces
and Complex Analysis 3, (1997), 152-156.

\bibitem{S}
Schwartz, L., Th\'eorie des distributions 1. Hermann, Paris, 1950.
-- 350p.

\bibitem{T1}
Tornehave, H., {\it Systems of zeros of holomorphic almost
periodic functions}, Kobenhavns Universitet Matematisk Institut,
Preprint no. 30, 1988.

\bibitem{T2}
Tornehave, H., {\it On the zeros of entire almost periodic
functions}, The Harald Bohr Centenary (Copenhagen 1987). Math.
Fys. Medd. Danske 42 (3), (1989), 125-142.

\end{thebibliography}
\end{document}